\documentclass[11pt]{amsart}
\usepackage[top=1.5in, bottom=1.5in, left=1.5in, right=1.5in]{geometry}
\expandafter\let\csname ver@amsthm.sty\endcsname\relax
\usepackage{amsmath}
\usepackage{amssymb}
\usepackage{hyperref}
\usepackage{amsthm}
\usepackage{mathtools}
\usepackage[capitalize,noabbrev]{cleveref}
\usepackage{xcolor}

\allowdisplaybreaks

\newtheorem{theorem}{Theorem}
\newtheorem{corollary}{Corollary}
\newtheorem{Example}{Example}
\newenvironment{example}
  {\begin{Example}\rm}{\end{Example}}
\newtheorem{Remark}{Remark}
\newenvironment{remark}
  {\begin{Remark}\rm}{\end{Remark}}
  \newtheorem{Question}{Question}
\newenvironment{question}
  {\begin{Question}\rm}{\end{Question}}
  
\crefname{Theorem}{Theorem}{Theorems}
\crefname{Corollary}{Corollary}{Corollaries}
\crefname{Example}{Example}{Examples}
\crefname{Remark}{Remark}{Remarks}
\crefname{Question}{Question}{Questions}

\newcommand{\emailhref}[1]{\email{\href{#1}{#1}}}
\newcommand{\dfn}[1]{\textcolor{blue}{\emph{#1}}}

\title{A note on M\"{o}bius functions of upho posets}
\author{Sam Hopkins}\emailhref{samuelfhopkins@gmail.com}
\address{Department of Mathematics, Howard University, Washington, DC 20059}

\begin{document}

\begin{abstract}
A poset is called upper homogeneous (or ``upho'') if every principal order filter of the poset is isomorphic to the whole poset. We observe that the rank and characteristic generating functions of upho posets are multiplicative inverses of one another.
\end{abstract}

\maketitle

We refer to~\cite[\S 3]{stanley2012ec1} for basic terminology and notation for posets. A poset $P$ is called \dfn{$\mathbb{N}$-graded} if we can write $P$ as a disjoint union $P= P_0\sqcup P_1 \sqcup P_2 \sqcup \cdots$ such that every maximal chain has the form $p_0 \lessdot p_1 \lessdot p_2 \lessdot \cdots$ with $p_i \in P_i$ for all $i$. The rank function $\rho\colon P\to \mathbb{N}$ of $P$ is then given by $\rho(p) = i$ if $p \in P_i$. We say $P$ has \dfn{finite type} if $\#P_i < \infty$ for all $i$. In this case we can form the \dfn{rank generating function}
\[ F_P(x) \coloneqq \sum_{i \geq 0} \#P_i \; x^i = \sum_{p\in P} x^{\rho(p)}.\]
Suppose further that $P$ has a minimum element $\hat{0}\in P_0$. Then we define
\[ \chi_P(x) \coloneqq \sum_{p\in P} \mu(\hat{0},p) \; x^{\rho(p)},\]
where $\mu(\cdot,\cdot)$ is the M\"{o}bius function of $P$. The analogous $\chi_P(x) \coloneqq  \sum_{p\in P} \mu(\hat{0},p) \; x^{\rho(p)}$ for a finite~$P$ (or, more often, its reciprocal polynomial) is called the characteristic polynomial of~$P$. So we refer to~$\chi_P(x)$ as the \dfn{characteristic generating function}.

The coefficients of these generating functions often have great combinatorial significance: for example, when $P=\Pi_n$ is the lattice of set partitions of $\{1,2,\ldots,n\}$ ordered by refinement, the coefficients of $\chi_P(x)$ and $F_P(x)$ are the Stirling numbers $s(n,k)$ and $S(n,k)$ of the 1st and 2nd kind, respectively.

In this note we observe that, for a special class of posets called ``upper homogeneous'' (or ``upho'' for short), there is a very simple relationship between the rank and characteristic generating functions: they are multiplicative inverses.

A poset $P$ is called \dfn{upper homogeneous (upho)} if we have $V_p \simeq P$ for all $p \in P$, where $V_p \coloneqq \{q\in P\colon q \geq p\}$ is the principal order filter (i.e., dual order ideal) generated by~$p$. Note that a nontrivial upho poset necessarily has a minimum, and is infinite if it has more than one element. Upho posets were introduced recently by Stanley~\cite{stanley2020upho} in his investigation of certain generating functions related to Stern's diatomic array~\cite{stanley2020stern} and the Fibonacci numbers~\cite{stanley2021rational}. Not much is known about the structure of upho posets in general, but see~\cite{gao2020upho} for a recent paper studying the rank generating functions of finite type $\mathbb{N}$-graded upho posets.

From now on, upho posets are assumed finite type $\mathbb{N}$-graded. Our main result is:

\begin{theorem} \label{thm:main}
For $P$ an upho poset, we have $F_P(x)=\chi_P(x)^{-1}$.
\end{theorem}

\begin{example} \label{ex:grid}
The ``grid'' $P=\mathbb{N}^n$ is upho with $F_P(x) = \frac{1}{(1-x)^n}$ and $\chi_P(x) = (1-x)^n$. (These computations follow immediately from the fact that for a Cartesian product $P_1\times P_2$ we have $F_{P_1\times P_2}(x)=F_{P_1}(x)F_{P_2}(x)$ and $\chi_{P_1\times P_2}(x) =\chi_{P_1}(x)\chi_{P_2}(x)$.)
\end{example}

\begin{example} \label{ex:tree}
The ``infinite (rooted) $n$-ary tree'' poset $P$ is upho with $F_P(x) = \frac{1}{1-nx}$ and $\chi_P(x) = 1-nx$.
\end{example}

\begin{example} \label{ex:bowtie}
Fix $n\geq 1$ and let $P$ be the $\mathbb{N}$-graded poset with $\#P_0=1$, $\#P_i=n$ for all $i\geq 1$, and all cover relations between any two adjacent ranks (the ``bowtie'' poset from \cite[Figure~1]{gao2020upho} is the case $n=2$ of this poset). Then $P$ is upho and has $F_P(x) = \frac{1+(n-1)x}{1-x}$ and $\chi_P(x)=\frac{1-x}{1+(n-1)x}$. 
\end{example}

\begin{proof}[Proof of \cref{thm:main}]
Let $P$ be upho. First we claim that for any $m\geq 0$,
\begin{equation} \label{eqn:chains}
 \sum_{i \geq 0} \#\{\textrm{chains $\hat{0}=p_0 < p_1 < \cdots < p_m$ of $P$} \colon \rho(p_m)=i \} \; x^i = (F_P(x)-1)^m
\end{equation}
Indeed, this is easily proved by induction: the number of ways to extend a chain $\hat{0}=p_0 < p_1 < \cdots < p_{m-1}$ with $\rho(p_{m-1}) = j$ to a chain $\hat{0}=p_0 < p_1 < \cdots < p_m$ with $\rho(p_{m}) = i$ is the coefficient of $x^{i-j}$ in $(F_P(x)-1)$, precisely because $V_{p_{m-1}}\simeq P$.

Next, we recall ``Philip Hall's theorem''~\cite[Proposition 3.8.5]{stanley2012ec1}, which says that a poset's M\"{o}bius function $\mu(\cdot,\cdot)$ satisfies
\[ \mu(p,q) = c_0 - c_1 + c_2 - c_3 + \cdots\]
where $c_i$ is the number of length $i$ chains $p=p_0 < p_1 < \cdots < p_i = q$ from $p$ to~$q$. 

Hence,
\begin{align*}
\chi_P(x) &= \sum_{p\in P} \mu(\hat{0},p) \; x^{\rho(p)} \\
&= \sum_{p\in P} \left( \sum_{m \geq 0} (-1)^m \#\{\textrm{chains $\hat{0}=p_0 < p_1 < \cdots < p_m=p$}\} \right) x^{\rho(p)} \\
&= \sum_{m \geq 0} (-1)^m  \sum_{i \geq 0} \#\{\textrm{chains $\hat{0}=p_0 < p_1 < \cdots < p_m$} \colon \rho(p_m)=i \} \; x^i \\
&= \sum_{m \geq 0} (-1)^m  (F_P(x)-1)^m = \frac{1}{1- (-(F_P(x)-1))} = F_P(x)^{-1},
\end{align*}
where from the 1st to the 2nd line we used Philip Hall's theorem, and from the 3rd to the 4th line we used~\eqref{eqn:chains}.
\end{proof}

\begin{remark} \label{rem:moebius}
An alternative proof of \cref{thm:main} is via M\"{o}bius inversion~\cite[\S 3.7]{stanley2012ec1}. Define $f(p) \coloneqq x^{\rho(p)}$ and $g(p) \coloneqq \sum_{q \geq p} f(q)$ for each $p\in P$. By M\"{o}bius inversion, $1=x^{\rho(\hat{0})} = \sum_{q \in P} \mu(\hat{0},q) \; g(q)$. But since $P$ is upho, $g(q) = x^{\rho(q)} \; F_P(x)$ for all $q \in P$, so that
\[ 1 = \sum_{q \in P} \mu(\hat{0},q) \; x^{\rho(q)} \; F_P(x) = F_P(x) \cdot \left( \sum_{q \in P} \mu(\hat{0},q) \; x^{\rho(q)}\right) = F_P(x) \cdot \chi_P(x).\]
In other words, $F_P(x) = \chi_P(x)^{-1}$. Because these sums are infinite, \cite[Proposition~3.7.2]{stanley2012ec1} as stated does not literally apply; nevertheless, these manipulations can still be justified by taking an appropriate limit in the ring of formal power series.
\end{remark}

M\"{o}bius functions are especially well behaved for lattices, so from now on we concentrate on the case of $P$ an upho lattice.

\begin{corollary} \label{cor:lattice}
Let $P$ be an upho lattice. Then $ F_P(x) = \chi_{P'}(x)^{-1}$, where 
\[P' \coloneqq \{p \in P\colon p \leq a_1 \vee a_2 \vee \cdots \vee a_k \textrm{ for some atoms $a_1,\ldots,a_k \in P$}\}\]
is the finite graded sub-lattice of elements below joins of atoms. The same is true if we replace ``lattice'' with ``meet semilattice'' everywhere.
\end{corollary}

\begin{proof}
By \cref{thm:main} it suffices to prove that $\chi_{P}(x) = \chi_{P'}(x)$. By supposition, an interval~$[\hat{0},p]$ for $p \in P$ is a finite lattice. Hence, by the crosscut theorem -- or specifically, its corollary~\cite[Corollary 3.9.5]{stanley2012ec1} -- we will have $\mu(\hat{0},p) = 0$ unless $p$ is a join of atoms. Thus, to record all non-zero M\"{o}bius function values we only need to consider intervals between~$\hat{0}$ and joins of atoms, so $\chi_{P}(x) = \chi_{P'}(x)$, as required. The only difference when~$P$ is a meet semilattice rather than a lattice is that some subsets of atoms may fail to have a join, but~$P'$ consists precisely of all elements below subsets of atoms which do have a join.
\end{proof}

\begin{remark}
A result of Gao, Guo, Seetharaman, and Seidel~\cite[Theorem 1.3]{gao2020upho} says that the rank generating function of a planar upho poset~$P$ (i.e., an upho poset whose Hasse diagram is planar) is the inverse of a polynomial. Every planar upho poset $P$ is a meet semilattice~\cite[Lemma 4.1]{gao2020upho}, so \cref{cor:lattice} is another way to see that its rank generating function is the inverse of a polynomial. In fact, Y.~Gao (private communication) pointed out that the M\"{o}bius function of a planar upho poset $P$ is
\[ \mu(\hat{0},p) = \begin{cases} 1 &\textrm{if $p=\hat{0}$ or $p$ is root-bifurcated,} \\
-1 &\textrm{if $p$ is an atom}, \\
0 &\textrm{otherwise}.\end{cases}\]
See~\cite[Definition~4.1]{gao2020upho} for the definition of root-bifurcated element, of which there are only finitely many~\cite[Lemma~4.5]{gao2020upho}.  Note also, by way of contrast, that Gao et al.~\cite[\S5]{gao2020upho} showed how rank generating functions of \emph{arbitrary} upho posets can be very complicated.
\end{remark}

\begin{remark} \label{rem:lattice}
\Cref{cor:lattice} says that a lot of information about the upho lattice $P$ is contained in the finite graded lattice $P'$ below the join of all atoms. But the whole structure of $P$ is not determined by $P'$. For example, as suggested in \cref{ex:grid}, with $P=\mathbb{N}^n$ we have $P'=$ the rank $n$ Boolean lattice. But a different $P$ with $P'=$ the rank $n$ Boolean lattice is given by $P=\{\textrm{finite } A\subseteq \{1,2,\ldots\}\colon \max(A) < \#A+n\}$ (with the order being inclusion). 
\end{remark}

In spite of the fact that the extension will not in general be unique, it is still natural to ask when one can ``go in the other direction'' and extend a $P'$ to a $P$.

\begin{question} \label{question:lattice}
Consider a finite graded lattice $P'$. Can one find an upho lattice $P$ such that $P'$ is the sub-lattice of $P$ below the join of all atoms?
\end{question}

\Cref{cor:lattice} says that for such a $P'$ to be extendable, it must be the case that $\chi_{P'}(x)^{-1}$ has all positive coefficients. So for a ``random'' $P'$ the answer to \cref{question:lattice} will be negative. On the other hand, in \cref{rem:lattice} we gave an affirmative answer when $P'=$ the rank $n$ Boolean lattice. We now review some other examples of well-studied finite graded lattices which give affirmative answers to \cref{question:lattice}.

\begin{example} \label{ex:prime}
Fix $n \geq 1$ and a prime $p$, and let~$P$ be the set of subgroups of $\mathbb{Z}^n$ of index a power of $p$ ordered by reverse inclusion. Then~$P$ is an upho lattice~\cite{stanley2020upho}, and $P' =$ the lattice of subspaces of $(\mathbb{Z}/p\mathbb{Z})^n$. One can compute directly (e.g.~using Hermite normal form) that $F_P(x)=\frac{1}{(1-x)(1-xp)\cdots (1-xp^{n-1})}$, or deduce this from the well-known formula $\chi_{P'}(x)=(1-x)(1-xp)\cdots(1-xp^{n-1})$ together with \cref{cor:lattice}.
\end{example}

\begin{example} \label{ex:partition}
Fix $n \geq 1$ and let $P$ be the poset whose elements are partitions of sets of the form~$\{1,2,\ldots,k\}$ (for some $k\geq n$) into $n$ blocks, with $\pi_1 \leq \pi_2$ if for every $B_1 \in \pi_1$ there is some $B_2 \in \pi_2$ with $B_1\subseteq B_2$. Then $P$ is an upho lattice (where the rank of a partition of $\{1,2,\ldots,k\}$ into $n$ blocks is $k-n$), and $P' = \Pi_{n+1}$. Again, one can compute directly $F_P(x) = \sum_{k \geq n} S(k,n)x^{k-n} = \frac{1}{(1-x)(1-2x)\cdots(1-nx)}$, or deduce this from $\chi_{P'}(x) = (1-x)(1-2x)\cdots(1-nx)$  together with \cref{cor:lattice}.
\end{example}

\begin{example} \label{ex:noncrossing}
V. Reiner (private communication) explained that taking $P$ to be the ``dual braid monoid'' of a finite Coxeter group~\cite{bessis2003dual}, we have $P' =$ the corresponding ``noncrossing partition lattice.'' The rank generating function and M\"{o}bius function connection for this particular example is explored in~\cite{josuatverges2021koszulality} (see also~\cite{cartier1969problemes}).
\end{example}

\begin{remark}
A finite graded poset $P$ of rank $n$ with a minimum and a maximum is called \dfn{uniform} if, for each $i=0,1,\ldots,n$, all principal order filters $V_p$ for $p \in P$ with $\rho(p)=n-i$ are isomorphic to the same fixed poset $Q_i$. The rank $n$ Boolean lattice, the lattice of subspaces of $(\mathbb{Z}/p\mathbb{Z})^n$, and $\Pi_{n+1}$ are all uniform. It is known (see~\cite[Theorem~6]{dowling1973geometric} and~\cite[Exercise 3.130(a)]{stanley2012ec1}) that, for such a $P$, the matrices of the 1st and 2nd kind Whitney numbers for these~$Q_i$ are inverses of one another. This generalizes the fact that the matrices $(s(i,j))^{i=1,\ldots,n}_{j=1,\ldots,n}$ and $(S(i,j))^{i=1,\ldots,n}_{j=1,\ldots,n}$ of the 1st and 2nd kind Stirling numbers are inverses. \Cref{thm:main} seems superficially quite similar to this fact about uniform posets, but we do not see any direct connection. However, it would definitely be reasonable to look at other sequences of uniform lattices in search of affirmative answers to \cref{question:lattice}. 
\end{remark}

We conclude with one additional, interesting corollary of \cref{thm:main}:

\begin{corollary} \label{cor:meet}
Let $P$ be an upho meet semilattice. Then for any~$m\geq 1$,
\[ \sum_{\substack{(p_1,\ldots,p_m)\in P^m, \\ p_1 \wedge \cdots \wedge p_m = \hat{0}}} x^{\rho(p_1)+\cdots+\rho(p_m)} = F_P(x)^m \cdot F_P(x^m)^{-1}. \]
\end{corollary}

\begin{proof}
This is what we get by combining~\cite[Exercise 3.89]{stanley2012ec1} and~\cref{thm:main}. Namely, for each $p\in P$, set
\[f(p) \coloneqq \sum_{\substack{(p_1,\ldots,p_m)\in P^m, \\ p_1 \wedge \cdots \wedge p_m = p}} x^{\rho(p_1)+\cdots+\rho(p_m)} \]
and $g(p) \coloneqq \sum_{q \geq p} f(q)$. Then by M\"{o}bius inversion
\[ \sum_{\substack{(p_1,\ldots,p_m)\in P^m, \\ p_1 \wedge \cdots \wedge p_m = \hat{0}}} x^{\rho(p_1)+\cdots+\rho(p_m)} =  f(\hat{0}) = \sum_{q \in P} \mu(\hat{0},q) \; g(q).\]
But since $P$ is upho, we have
\[g(q) =\sum_{\substack{(p_1,\ldots,p_m)\in P^m, \\ p_1 \wedge \cdots \wedge p_m \geq q}} x^{\rho(p_1)+\cdots+\rho(p_m)} = \sum_{\substack{(p_1,\ldots,p_m)\in P^m, \\ p_1, \ldots, p_m \geq q}} x^{\rho(p_1)+\cdots+\rho(p_m)} =  (x^{\rho(q)} \; F_P(x))^m\]
for all $q\in P$, so that
\begin{align*}
\sum_{\substack{(p_1,\ldots,p_m)\in P^m, \\ p_1 \wedge \cdots \wedge p_m = \hat{0}}} x^{\rho(p_1)+\cdots+\rho(p_m)} &=  \sum_{q \in P} \mu(\hat{0},q)  \; (x^{\rho(q)} \; F_P(x))^m \\
&= F_P(x)^m \cdot \chi_P(x^m) = F_P(x)^m \cdot F_P(x^m)^{-1},
\end{align*} 
where in the last line we applied~\cref{thm:main}.
\end{proof}

\Cref{cor:meet} could in theory be useful for addressing \cref{question:lattice}. As mentioned, already \cref{cor:lattice} implies that for a finite graded lattice $P'$ to be extendable to an upho lattice $P$,  $\chi_{P'}(x)^{-1}$ must have all positive coefficients. \Cref{cor:meet} says that additionally $\chi_{P'}(x)^{-m} \cdot \chi_{P'}(x^m)$ must have all positive coefficients, for all $m\geq 1$.

\subsection*{Acknowledgments} I thank Yibo Gao, Vic Reiner, and Richard Stanley for helpful comments, and the anonymous referee for comments that improved the exposition.

\bibliography{upho_moebius}{}
\bibliographystyle{abbrv}

\end{document}